\documentclass[11pt]{article}

\usepackage{fullpage} % Package to use full page
\usepackage{parskip} % Package to tweak paragraph skipping
\usepackage{tikz} % Package for drawing

\usepackage{amsmath}
\usepackage{amsthm}

\usepackage{hyperref}
\usepackage{listings}
\usepackage{cite,graphics}

\usepackage{algorithm,algpseudocode}
\usepackage{amssymb}

\usepackage{subfigure}
\usepackage{diagbox}
\usepackage{tikz}
\usetikzlibrary{positioning}

\tikzstyle{inputNode}=[draw,circle,minimum size=12pt,inner sep=0pt]
\tikzstyle{stateTransition}=[-stealth, thick]

  \def\and{\  and }

\title{A Local Deep Learning Method for Solving High Order Partial Differential Equations}

\author{Jiang Yang\thanks{SUSTech International Center for Mathematics \& Department of Mathematics, Southern University of Science and Technology, and Guangdong Provincial Key Laboratory of Computational Science and Material Design, Southern University of Science and Technology, Shenzhen, China, {\tt yangj7@sustech.edu.cn}}
\and
Quanhui Zhu\thanks{Department of Mathematics, Southern University of Science and Technology, Shenzhen, China
{\tt 11849393@mail.sustech.edu.cn}}
}

\begin{document}

\maketitle
\begin{abstract}
    At present, deep learning based methods are being employed to resolve the computational challenges of high-dimensional partial differential equations (PDEs). But the computation of the high order derivatives of neural networks is costly, and high order derivatives lack robustness for training purposes. We propose a novel approach to solving PDEs with high order derivatives by simultaneously approximating the function value and derivatives. We introduce intermediate variables to rewrite the PDEs into a system of low order differential equations as what is done in the local discontinuous Galerkin method. The intermediate variables and the solutions to the PDEs are simultaneously approximated by a multi-output deep neural network. By taking the residual of the system as a loss function, we can optimize the network parameters to approximate the solution. The whole process relies on low order derivatives. Numerous numerical examples are carried out to demonstrate that our local deep learning is efficient, robust, flexible, and is particularly well-suited for high-dimensional PDEs with high order derivatives.
\end{abstract}


\section{Introduction}
Partial differential equations (PDEs) play a significant role in the fields of physics, chemistry, biology, engineering, finance, and others. Classical numerical methods focus on designing efficient, accurate, and stable numerical schemes. Within the context of high-dimensional problems, however, the curse of dimensionality renders classical numerical methods impractical. As a result, many mathematicians have introduced neural networks into PDEs precisely because multilayer feedforward networks are proven to be universal approximators for the PDEs\cite{Hornik1989Multilayer,Hornik1991Approximation}. More specifically, once the network structure is determined, any order derivatives of the neural network can be obtained analytically. Coupled with the automatic differentiation technique, neural networks can be applied to solve PDEs\cite{DBLP:journals/corr/BaydinPR15}. Depending upon different purposes, neural networks can be used to approximate the solution function, represent the solution solver, and even invert the equations.

In this paper, we consider the deep learning method as a means to solve the following $k$-th order initial boundary value problem (IBVP)
\begin{equation}
    \left\{
    \begin{array}{rcll}
        u_t &=& \mathcal{L}(u), &x\in\Omega,t\in[0,T],\\
        u(x,0) &=& u_0(x), &x\in\Omega,\\
        \mathcal{B}u &=& \mathbf{g}, & x\in\partial\Omega,t\in[0,T],
    \end{array}
    \right.
    \label{PDE}
\end{equation}
where $\Omega\subset\mathbb{R}^d$, $d\in \mathbb{N}_+$, $\mathcal{L}(u) = F(x,t,u,Du , \cdots, D^ku)$, $F$ and $\mathbf{g}$ are linear or nonlinear functions, $\mathcal{B}$ is the boundary condition operator, and the $p$-th order derivative operator $D^p$ consists of
\begin{equation*}
   \partial _{x_1}^{\alpha_1}\partial _{x_2}^{\alpha_2}\cdots\partial _{x_d}^{\alpha_d}u,\;\text{with}\; \sum \alpha_i=p,\;\alpha_i\in \mathbb{N}.
\end{equation*}
The neural network function $\varphi(x,t;\theta):\mathbb{R}^{d+1} \times \Theta^M \mapsto \mathbb{R}^m$, is defined as follows
\begin{equation}
    \begin{array}{rcll}
        \varphi(x,t;\theta) &=& \mathcal{N}_{\text{out}} \circ \mathcal{N}_L \circ \cdots \circ \mathcal{N}_1 \circ \mathcal{N}_{\text{in}}(x,t),\\
        \mathcal{N}_{\text{in}}(x,t) &=& \sigma_{\text{in}}(\alpha x +\beta t + b), \alpha\in\mathbb{R}^{n\times d},\beta,b\in\mathbb{R}^{n}, \\
        \mathcal{N}_{\text{out}}(y) &=& \sigma_{\text{out}}(\gamma y + c), \gamma\in\mathbb{R}^{m\times n}, c\in\mathbb{R}^m,
    \end{array}
    \label{fr::L-DNN}
\end{equation}
where $d$ is the dimension of $x$, $m$ is the dimension of the output, $n$ is the width of the hidden layers, $\sigma$ is a nonpolynomial active function, $L$ is the number of the hidden layers (i.e., the network's depth) and $\mathcal{N}_i:\mathbb{R}^n\mapsto\mathbb{R}^n$ is the structure of the hidden layers. Our goal is to find a suitable neural network $\varphi(x,t;\theta)$ to approximate a solution $u(x,t)$ to the problem (\ref{PDE}).

In most of existing literatures, the loss function is determined by either the PDEs or an equivalent formulation. For instance, the parabolic PDE is reformulated as a backward stochastic differential equation in \cite{Han2018Solving,weinan2017deep,han2020derivative}, where the loss function is given by the solution of the backward stochastic differential equation, and the training process is shown to be a deep reinforcement learning process. In \cite{Sirignano2017DGM,RAISSI2019686}, the solution is approximated by a neural network. The proposal of a mesh-free algorithm makes high-dimensional calculations feasible. \cite{ZANG2020109409} considers variational problems, and the loss function is defined as a weak formulation. Further, \cite{DBLP:journals/corr/abs-1710-00211} introduces an adversarial network as a test function in variational problems; this is particularly suitable for high-dimensional PDEs defined in irregular domains.

We consider using deep learning methods to solve PDEs with high order derivatives. The cost of computing high order derivatives for neural networks is prohibitive. Addressing the issue of deep learning methods, \cite{Sirignano2017DGM} proposes a Monte Carlo method to approximate second order derivatives. In \cite{DBLP:journals/corr/abs-1710-00211,ZANG2020109409}, the variational form reduces the order of derivatives through the integration by parts. \cite{han2020derivative} proposes a derivative free method for parabolic PDEs by solving the equivalent BSDE problem. But, for the higher order derivatives of neural networks, there are still numerous computational challenges. High order derivatives limits the choices of network structure, influences the robustness in training, and are expansive to be computed.

The local discontinuous Galerkin (LDG) method introduces new variables and rewrites the problem (\ref{PDE}) as a system of first order differential equations \cite{cockburn1998local,yan2002local,yan2002ldg,xu2010local}. Then, the method is obtained by discretizing the system with the discontinuous Galerkin method. The reduction of order technique in LDG  inspires us to use a similar technique to compute high order derivatives in deep learning. To this end,
we first rewrite the PDEs to a system of low order differential equations. A neural network with multiple outputs is then used to approximate the solution and intermediate variables. Taking the $L_2$ residual of the system as the loss function, we can optimize the neural network to approximate the solution of (\ref{PDE}). Unlike the classical deep learning methods, our approach avoids calculating the high order derivatives. As consequence, it is more efficient and stable.

Our paper is organized as follows. In Section \ref{sec2}, we briefly introduce the deep learning method for solving PDEs and illustrate the difficulties in computing the high order derivatives of the neural networks. After rewriting the problem as a system of low order equations, the local deep Galerkin method and the local deep Ritz method are proposed in Section \ref{sec3}. The advantages of this method are provided in Section \ref{sec4}, and the corresponding numerical experiments are presented in Section \ref{sec5}. Several concluding remarks are given in the final section.

\section{Preliminaries}\label{sec2}

In this section, we present a deep learning based method for solving PDEs and show the main difficulties in computing the high order derivatives of neural networks.

\subsection{Deep Learning Based Method for Solving PDEs}
A deep neural network defined as (\ref{fr::L-DNN}) is used to approximate the solution of (\ref{PDE}). Substituting the neural network into (\ref{PDE}), we have
\begin{equation}
    \left\{
    \begin{array}{rcll}
        \varphi_t &=& \mathcal{L}(\varphi), &x\in\Omega,t\in[0,T],\\
        \varphi(x,0) &=& u_0(x), &x\in\Omega,\\
        \mathcal{B}\varphi&=& \mathbf{g}, &x\in\partial\Omega,t\in[0,T].
    \end{array}
    \right.
    \label{neural PDE}
\end{equation}
Instead of solving the equation step by step under a given initial value, the neural network parameters should be optimized to satisfy the dynamic system, the initial value condition and the boundary condition. The loss function $J(\varphi)$ is defined by the $L_2$ norm mostly\cite{Sirignano2017DGM,RAISSI2019686,hayati2007feedforward,lagaris1998artificial,dockhorn2019discussion,lagaris2000neural,lu2019deepxde}, i.e.,
\begin{equation}
    J(\varphi) = \|\varphi_t - \mathcal{L}(\varphi)\|^2_{\Omega_T} + \|\varphi(\cdot,0) - u_0(\cdot)\|^2_{\Omega} + \|\mathcal{B}\varphi-\mathbf{g}\|^2_{\partial\Omega_T},
\end{equation}
where $\Omega_T = \Omega\times[0,T]$, $\partial\Omega_T = \partial\Omega\times[0,T]$.
Denote $J(\varphi) := J_e(\varphi) + J_i(\varphi) + J_b(\varphi)$, which represents the equation loss, the initial condition loss, and the boundary condition loss, respectively. The learning algorithm can be described as follows:
\begin{itemize}
    \item [1.] Build up a neural network $\varphi(x,t;\theta)$, which determines active functions, the hidden layer structure, and the network's depth and width. The inputs are space $x$ and time $t$ and the output is the function value at $(x,t)$, which approximates the solution of (\ref{PDE}). $\theta$ is the trainable parameters in the neural network.
    \item [2.] Obtain random samples $\mathcal{D}$ from within the domain $\Omega_T$, $\partial\Omega_T$ and $\{0\}\times\Omega$. Random sampling makes high-dimensional calculations feasible. The mesh-free property is one of the most critical differences between deep learning methods and classical numerical methods.
    \item [3.] Solve the optimization problem on the given sampling set $\mathcal{D}$:
    \begin{equation*}
        \min_\theta J(\varphi)|_{\mathcal{D}}.
    \end{equation*}
  The discrete form is given as
    \begin{equation}
        J(\varphi)|_{\mathcal{D}} = \frac{w_e}{N_e}\sum_{(x_{e},t_{e})}(\varphi_t-\mathcal{L}(\varphi))^2 + \frac{w_i}{N_i}\sum_{(x_{i},t_{i})}(\varphi-u_0)^2 +\frac{w_b}{N_b}\sum_{(x_{b},t_{b})}\|\mathcal{B}\varphi-\mathbf{g}\|_2^2,
        \label{eq::general_loss}
    \end{equation}
    where $w$ is the weight of per term and $N$ is the number of sampling nodes.
     Finding the optimal parameters for a fixed width neural network is difficult since the optimization problem is nonconvex. The Adam optimizer is a popular choice for deep learning\cite{kingma2014adam}.
    \item [4.] Repeat steps 2 and 3 until the result converges.
\end{itemize}
The selection of the active function, initialization method, random sampling distribution, and optimization method will affect the approximation of the network. Incorrect settings will result in either the failure of the the neural network to converge or a very slow rate of convergence.

\subsection{High Order Derivatives of the Deep Neural Network in PDEs}
\label{sec::complexity}
 High order derivatives of the neural network rarely appear in classic deep learning problems. But, to solve high order PDEs, we have to calculate the value of high order derivatives. Consider a fully connected neural network $\varphi(x,\theta):\mathbb{R}^d\times\Theta^M\mapsto\mathbb{R}$ similar to Fig. (\ref{st::fc}), in which there are $L$ hidden layers and $n$ neurons per layer. The layer  structure is
\begin{equation}
    \mathcal{N}_i(x) = \sigma(W_i\mathcal{N}_{i-1}(x) + b_i), W_i\in\mathbb{R}^{n\times n}, b_i\in\mathbb{R}, i=2\cdots,L.
\end{equation}
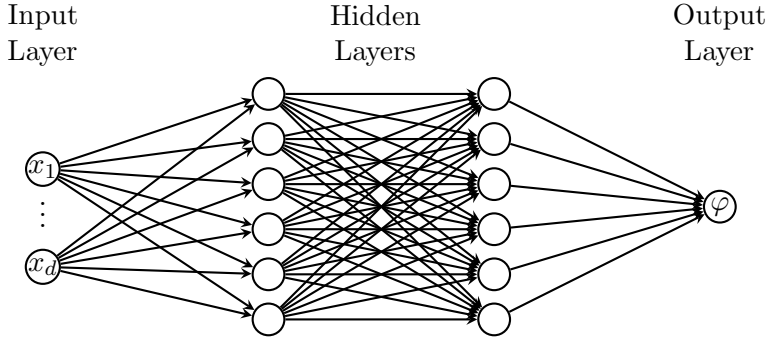
\begin{figure}[H]
    \begin{center}
    \begin{tikzpicture}
        \node[inputNode, thick] (x1) at (1, 0.5) {$x_1$};
        \node[inputNode, thick] (x2) at (1, -0.8) {$x_d$};
        \node[] (x3) at (1, 0) {$\vdots$};

        \node[inputNode, thick] (h5) at (4, 1.5) {};
        \node[inputNode, thick] (h6) at (4, -1.5) {};
        \node[inputNode, thick] (h1) at (4, 0.9) {};
        \node[inputNode, thick] (h2) at (4, 0.3) {};
        \node[inputNode, thick] (h3) at (4, -0.3) {};
        \node[inputNode, thick] (h4) at (4, -0.9) {};

        \node[inputNode, thick] (d1) at (7, 1.5) {};
        \node[inputNode, thick] (d2) at (7, 0.9) {};
        \node[inputNode, thick] (d3) at (7, 0.3) {};
        \node[inputNode, thick] (d4) at (7, -0.3) {};
        \node[inputNode, thick] (d5) at (7, -0.9) {};
        \node[inputNode, thick] (d6) at (7, -1.5) {};

        \node[inputNode, thick] (o1) at (10, 0) {$\varphi$};

        \draw[stateTransition] (x1) -- (h1);
        \draw[stateTransition] (x1) -- (h2);
        \draw[stateTransition] (x1) -- (h3);
        \draw[stateTransition] (x1) -- (h4);
        \draw[stateTransition] (x1) -- (h5);
        \draw[stateTransition] (x1) -- (h6);
        \draw[stateTransition] (x2) -- (h1);
        \draw[stateTransition] (x2) -- (h2);
        \draw[stateTransition] (x2) -- (h3);
        \draw[stateTransition] (x2) -- (h4);
        \draw[stateTransition] (x2) -- (h5);
        \draw[stateTransition] (x2) -- (h6);

        \draw[stateTransition] (h1) -- (d1);
        \draw[stateTransition] (h1) -- (d2);
        \draw[stateTransition] (h1) -- (d3);
        \draw[stateTransition] (h1) -- (d4);
        \draw[stateTransition] (h1) -- (d5);
        \draw[stateTransition] (h1) -- (d6);
        \draw[stateTransition] (h2) -- (d1);
        \draw[stateTransition] (h2) -- (d2);
        \draw[stateTransition] (h2) -- (d3);
        \draw[stateTransition] (h2) -- (d4);
        \draw[stateTransition] (h2) -- (d5);
        \draw[stateTransition] (h2) -- (d6);
        \draw[stateTransition] (h3) -- (d1);
        \draw[stateTransition] (h3) -- (d2);
        \draw[stateTransition] (h3) -- (d3);
        \draw[stateTransition] (h3) -- (d4);
        \draw[stateTransition] (h3) -- (d5);
        \draw[stateTransition] (h3) -- (d6);
        \draw[stateTransition] (h4) -- (d1);
        \draw[stateTransition] (h4) -- (d2);
        \draw[stateTransition] (h4) -- (d3);
        \draw[stateTransition] (h4) -- (d4);
        \draw[stateTransition] (h4) -- (d5);
        \draw[stateTransition] (h4) -- (d6);
        \draw[stateTransition] (h5) -- (d1);
        \draw[stateTransition] (h5) -- (d2);
        \draw[stateTransition] (h5) -- (d3);
        \draw[stateTransition] (h5) -- (d4);
        \draw[stateTransition] (h5) -- (d5);
        \draw[stateTransition] (h5) -- (d6);
        \draw[stateTransition] (h6) -- (d1);
        \draw[stateTransition] (h6) -- (d2);
        \draw[stateTransition] (h6) -- (d3);
        \draw[stateTransition] (h6) -- (d4);
        \draw[stateTransition] (h6) -- (d5);
        \draw[stateTransition] (h6) -- (d6);

        \draw[stateTransition] (d1) -- (o1);
        \draw[stateTransition] (d2) -- (o1);
        \draw[stateTransition] (d3) -- (o1);
        \draw[stateTransition] (d4) -- (o1);
        \draw[stateTransition] (d5) -- (o1);
        \draw[stateTransition] (d6) -- (o1);

        \node[above=of x1, align=center] (l1) {Input \\ Layer};
        \node[right=8em of l1, align=center] (l2) {Hidden \\ Layers};
        \node[right=8em of l2, align=center] (l3) {Output \\ Layer};

    \end{tikzpicture}
\end{center}
\caption{The structure of the fully connected neural network.}
\label{st::fc}
\end{figure}
Based on the chain rule, the computational cost of $k$-th order derivative of $\varphi$ is about $O(L^kn^{2k})$, and $D_x^{k} \varphi$ costs about $O(L^kd^kn^{2k})$. Although the automatic differentiation technique provides some convenience in practical problems, the exponential growth of the order $k$ is unacceptable for solving a specific high order PDE.

Putting aside the problem of computational efficiency, there are still inherent challenges to using deep neural networks to solve PDEs with high order derivatives. The gradient vanishing and exploding problems of neural networks have been discussed for many years\cite{bengio1994learning,hochreiter1998vanishing,pascanu2013difficulty,hanin2018neural,grosse2017lecture}. And each probably makes an appearance when PDEs are being solved. For simplicity, we consider when $n=d=1$ and $\mathcal{N}_{\text{out}}(x)=\mathcal{N}_{\text{in}}(x)=x$, i.e.,
\begin{equation}
    \varphi(x) = \prod_{i=1}^L \sigma(\mathcal{N}_i).
    \label{eq::onelayer}
\end{equation}
Then, we have the following first order derivative
\begin{equation}
    \frac{\partial \varphi}{\partial x} = \frac{\partial \varphi}{\partial \mathcal{N}_L}\cdot \prod_{i=2}^L \frac{\partial \mathcal{N}_i}{\partial \mathcal{N}_{i-1}}
    \cdot \frac{\partial \mathcal{N}_{1}}{\partial x}= \prod_{i=1}^L W_i\sigma'(\mathcal{N}_i).
\end{equation}
Given the active function $\sigma(x) = \tanh(x)$, we have $\sigma'(x)=1-\tanh^2(x) < 1$, except the point $x=0$. With a normal initialization $|W|<1$, $\varphi_x(x;\theta)\sim W^L\sigma'(x)^L$ will become small and this leads to the gradient vanishing problem. Similarly, if the active function satisfies $|W\sigma'(x)|>1$, it will result in the gradient exploding problem.

There are many mature deep learning techniques for resolving these problems in classification, regression, and so on. But, for deep learning methods in PDEs, things are different. Consider a $k$-th order ordinary differential equation
\begin{equation}
    F(x,u, u',u'',\cdots, u^{(k)})=0.
    \label{eq::ODE}
\end{equation}
We use the above neural network to approximate the solution $u(x)$. Substituting (\ref{eq::onelayer}) into (\ref{eq::ODE}) and denoting $\sigma^{(k)}_i = \sigma^{(k)}(\mathcal{N}_i), k\in \mathbb{N}$, we have
\begin{equation}
    F(x, \prod_{i=1}^L \sigma_i, \prod_{i=1}^L W_i\sigma'_i,\cdots, \prod_{i=1}^L W^k_i\sigma^{(k)}_i)=0.
\end{equation}

\begin{figure}[htbp]
    \centering
    \includegraphics[width=0.8\textwidth]{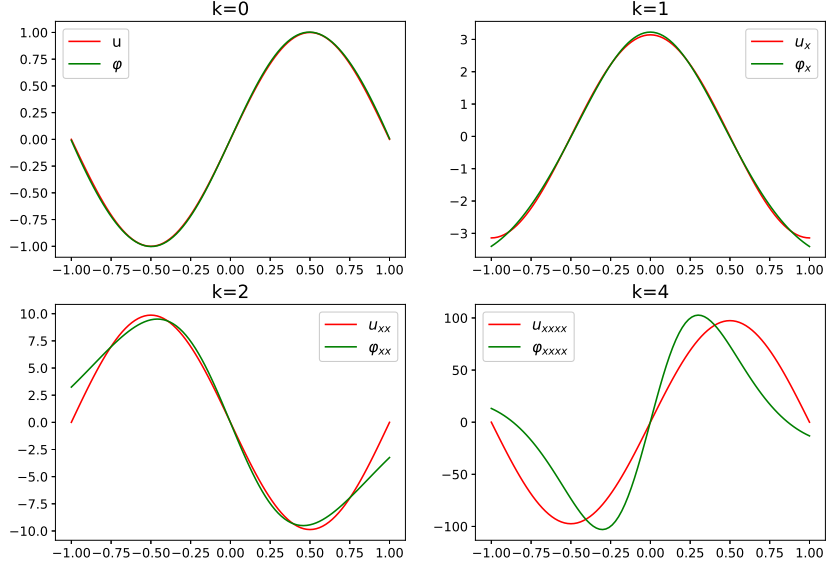}\vspace{-1em}
    \caption{Different order derivatives of the neural network $D^k \varphi$ compared to the target function $D^k u$ in one dimension.}
    \label{fig::gradient_scale}
\end{figure}

The issue is different from that of the classical deep learning problems. Our concern cannot just be limited to whether $|W\sigma'|$ is exploding or vanishing. The high order derivatives of the active function and the high order powers of parameters bring new difficulties. Taking $y:=\sigma(x)=\tanh(x)$ as an example, we have
\begin{equation}
    \begin{aligned}
        y' &= 1 - y^2, &y'\in&(0,1],\\
        y''&= -2y(1-y^2),&y''\in&[-\frac{4\sqrt{3}}{9},\frac{4\sqrt{3}}{9}],\\
        y^{(3)}&= (6y^2-2)(1-y^2),&y^{(3)}\in&[-2,\frac{2}{3}],\\
        y^{(4)}&= 8y(2-3y^2)(1-y^2),&y^{(4)}\in&(-4.086,4.086).
    \end{aligned}
\end{equation}
The derivative value can be greater than $1$ and the different order derivatives are controlled by the different scales. Fig. \ref{fig::gradient_scale} gives an example of gradient exploding problems of high order derivatives in using a $L=3,n=32$ neural network $\varphi$ to approximate the function $u(x)=\sin(\pi x)$. When order $k=4$, the gap between two derivatives is about $O(10)$, even though $\varphi(x)$ appproximates $u(x)$ well.

We use the following terms to represent the different order derivatives of a neural network with $L$ hidden layers
\begin{equation}
    \varphi^{(k)}(x;\theta) \sim (\sigma^{(k)}(x)\theta^k)^L.
\end{equation}
The gradient of the parameters is given as
\begin{equation}
    \nabla_\theta \varphi^{(k)}(x;\theta) \sim (x\sigma^{(k+1)}(x)\theta^k+k\sigma^{(k)}(x)\theta^{k-1})^L.
\end{equation}
The restriction on $\sigma'(x)\theta$ cannot restrict higher order derivatives. High order derivatives become more sensitive and unstable. It seems that different order derivatives are learned on different scales. When different order derivatives are in same equation, optimizing the loss is difficult (, and see also Section \ref{sec::CH}). The neural networks lack the efficiency and robustness of computing high order derivatives.

\section{Methodology}\label{sec3}
In this section, we propose employing local deep learning methods (LDLM) to overcome the derivative calculation problem of neural networks by using a multi-output neural network and a loss function of the equivalent system.

\subsection{The System of PDEs}
Consider a $k$-th order IBVP
\begin{equation}
    \left\{
    \begin{array}{rcll}
        u_t &=& F(u,D u, D^2 u, \cdots, D^k u), &x \in \Omega, t\in [0,T],\\
        u(x,0) &=& u_0(x), &x\in \Omega, \\
        \mathcal{B}u &=& \mathbf{g}, &x\in \partial\Omega,t\in[0,T].
    \end{array}
    \right.
    \label{system}
\end{equation}
Similar to the first part of the local discontinuous Galerkin method, we introduce the intermediate variables $\{v_i\}_{i=1}^{k}$, where $v_i\in \mathbb{R}^{d^i}, i=1,\cdots,k$. Then the PDE can be rewritten as the following system
\begin{equation}
    \left\{
        \begin{array}{rcll}
            u_t &=& F(u,v_1,v_2,\cdots,v_{k}), &x\in\Omega,t\in[0,T],\\
            v_1 &=& D u, & x\in\Omega,t\in[0,T],\\
            v_{i+1} &=& D v_i, & i=1,\cdots,k-1,x\in\Omega,t\in[0,T],\\
            u(x,0)& = & u_0(x), &x\in\Omega,\\
            \mathcal{B}u&= & \mathbf{g},&x\in\partial \Omega,t\in[0,T].
        \end{array}
    \right.
\end{equation}
Only the first order derivatives are included in the system. Similarly, we can build a system of the second order PDEs.
\begin{equation}
    \left\{
        \begin{array}{rcll}
            u_t &=& F(u,Du,\cdots,w_{[\frac{k}{2}]},Dw_{[\frac{k}{2}]}), &x\in\Omega,t\in[0,T],\\
            w_1 &=& D^2u, & x\in\Omega,t\in[0,T],\\
            w_{i+1} &=& D^2 w_i, & i=1,\cdots,[\frac{k}{2}],x\in\Omega,t\in[0,T],\\
            u(x,0)& = & u_0(x), &x\in\Omega,\\
            \mathcal{B}u&= & \mathbf{g},&x\in\partial \Omega,t\in[0,T].
        \end{array}
    \right.
    \label{system2}
\end{equation}
Some specific examples are given in Table \ref{tb::system}.
\begin{table}[htbp]
    \centering
	\caption{The System Form of Serval Classical Equations}
    \begin{tabular}{llll}
    \hline
    Equation& Origin Form & First Order System & Second Order System \\
    \hline
    Heat & $u_t=\Delta u$ & $\left\{\begin{aligned}
        u_t &= \nabla \cdot v\\ v &= \nabla u
    \end{aligned}\right.$ & $\left\{\begin{aligned}
        u_t &= v\\ v &= \Delta u
    \end{aligned}\right.$\\
    \hline
    Allen--Cahn& $u_t=\epsilon \Delta u + f(u)$& $\left\{\begin{aligned}
        u_t &= \epsilon \nabla\cdot v + f(u)\\ v &= \nabla u
    \end{aligned}\right.$ & $\left\{\begin{aligned}
        u_t &= \epsilon v + f(u)\\v&=\Delta u
    \end{aligned}\right.$\\
    \hline
    Cahn--Hilliard& $u_t = -\Delta(\epsilon\Delta u + f(u))$ & $\left\{\begin{aligned}
        u_t &= -\nabla \cdot v\\ v &= \nabla \phi \\ \phi&=\epsilon \nabla\cdot w + f(u)\\ w &= \nabla u
    \end{aligned}\right.$ &$\left\{\begin{aligned}
        u_t &= -\Delta v\\ v&= \epsilon\Delta u + f(u)
    \end{aligned}\right. $\\
    \hline
    KdV&$u_t + 6uu_x +u_{xxx} =0$&$\left\{\begin{aligned}
        u_t+6uv+w_x&=0\\v-u_x&=0 \\w-v_x&=0
    \end{aligned}\right.$&$\left\{\begin{aligned}
        u_t+6uu_x+v_x&=0\\v-u_{xx}&=0
    \end{aligned}\right. $
    \\
    \hline
    \end{tabular}
    \label{tb::system}
\end{table}

It is intuitive that the introduction of intermediate variables can effectively reduce the order of the derivatives.

\subsection{Multi-Output Neural Network}
\begin{figure}[H]

\begin{center}
    \begin{tikzpicture}
        \node[inputNode, thick] (x1) at (1, 0.5) {$x$};
        \node[inputNode, thick] (x2) at (1, -0.5) {$t$};

        \node[inputNode, thick] (h1) at (4, 0.9) {};
        \node[inputNode, thick] (h2) at (4, 0.3) {};
        \node[inputNode, thick] (h3) at (4, -0.3) {};
        \node[inputNode, thick] (h4) at (4, -0.9) {};

        \node[inputNode, thick] (d1) at (7, 1.5) {};
        \node[inputNode, thick] (d2) at (7, 0.9) {};
        \node[inputNode, thick] (d3) at (7, 0.3) {};
        \node[inputNode, thick] (d4) at (7, -0.3) {};
        \node[inputNode, thick] (d5) at (7, -0.9) {};
        \node[inputNode, thick] (d6) at (7, -1.5) {};

        \node[inputNode, thick] (o1) at (10, 0.75) {$u$};
        \node[inputNode, thick] (o2) at (10, 0) {$v_1$};
        \node[inputNode, thick] (o3) at (10, -0.75) {$v_2$};

        \draw[stateTransition] (x1) -- (h1);
        \draw[stateTransition] (x1) -- (h2);
        \draw[stateTransition] (x1) -- (h3);
        \draw[stateTransition] (x1) -- (h4);
        \draw[stateTransition] (x2) -- (h1);
        \draw[stateTransition] (x2) -- (h2);
        \draw[stateTransition] (x2) -- (h3);
        \draw[stateTransition] (x2) -- (h4);

        \draw[stateTransition] (h1) -- (d1);
        \draw[stateTransition] (h1) -- (d2);
        \draw[stateTransition] (h1) -- (d3);
        \draw[stateTransition] (h1) -- (d4);
        \draw[stateTransition] (h1) -- (d5);
        \draw[stateTransition] (h1) -- (d6);
        \draw[stateTransition] (h2) -- (d1);
        \draw[stateTransition] (h2) -- (d2);
        \draw[stateTransition] (h2) -- (d3);
        \draw[stateTransition] (h2) -- (d4);
        \draw[stateTransition] (h2) -- (d5);
        \draw[stateTransition] (h2) -- (d6);
        \draw[stateTransition] (h3) -- (d1);
        \draw[stateTransition] (h3) -- (d2);
        \draw[stateTransition] (h3) -- (d3);
        \draw[stateTransition] (h3) -- (d4);
        \draw[stateTransition] (h3) -- (d5);
        \draw[stateTransition] (h3) -- (d6);
        \draw[stateTransition] (h4) -- (d1);
        \draw[stateTransition] (h4) -- (d2);
        \draw[stateTransition] (h4) -- (d3);
        \draw[stateTransition] (h4) -- (d4);
        \draw[stateTransition] (h4) -- (d5);
        \draw[stateTransition] (h4) -- (d6);

        \draw[stateTransition] (d1) -- (o1);
        \draw[stateTransition] (d2) -- (o1);
        \draw[dashed] (d3) -- (o1);
        \draw[dashed] (d4) -- (o1);
        \draw[dashed] (d5) -- (o1);
        \draw[dashed] (d6) -- (o1);
        \draw[dashed] (d1) -- (o2);
        \draw[dashed] (d2) -- (o2);
        \draw[stateTransition] (d3) -- (o2);
        \draw[stateTransition] (d4) -- (o2);
        \draw[dashed] (d5) -- (o2);
        \draw[dashed] (d6) -- (o2);
        \draw[dashed] (d1) -- (o3);
        \draw[dashed] (d2) -- (o3);
        \draw[dashed] (d3) -- (o3);
        \draw[dashed] (d4) -- (o3);
        \draw[stateTransition] (d5) -- (o3);
        \draw[stateTransition] (d6) -- (o3);

        \node[above=of x1, align=center] (l1) {Input \\ Layer};
        \node[right=8em of l1, align=center] (l2) {Hidden \\ Layers};
        \node[right=8em of l2, align=center] (l3) {Output \\ Layer};

    \end{tikzpicture}
\end{center}
\caption{The multi-output neural network for approximating all of the intermediate variables.}
\label{st::1d}
\end{figure}
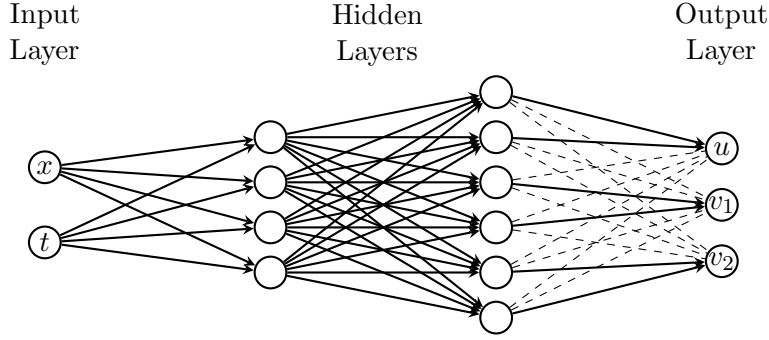

For approximating the intermediate variables and solution, a multi-output neural network $\varphi(x,t;\theta):\mathbb{R}^{d+1}\times \Theta^M\mapsto \mathbb{R}^m $ is needed. The 1-D coupled neural network is described in Fig. \ref{st::1d}.

In constrast to previous methods, all necessary intermediate variables are included in the output. For greater accuracy, the output layer can be a series of hidden layers. The final active function of each output is usually uniquely determined by the problem. When the number of intermediate variables increases, we only need to change the width of the output layer, and this is much cheaper than computing derivatives.

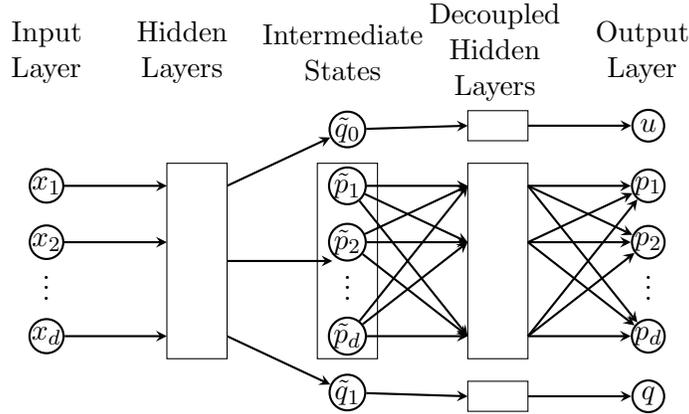
\begin{figure}[H]
\begin{center}
    \begin{tikzpicture}
        \node[inputNode, thick] (x1) at (0, 1) {$x_1$};
        \node[inputNode, thick] (x2) at (0, 0.25) {$x_2$};
        \node (dots) at (0,-0.25) {$\vdots$};
        \node[inputNode, thick] (x4) at (0, -1) {$x_d$};

        \draw[stateTransition] (x1) -- (1.6,1);
        \draw[stateTransition] (x2) -- (1.6,0.25);
        \draw[stateTransition] (x4) -- (1.6,-1);

        \draw (1.6,1.3) rectangle (2.4,-1.3);

        \node[inputNode, thick] (p1) at (4, 1.75) {$\tilde{q}_0$};
        \node[inputNode, thick] (p2) at (4, 1) {$\tilde{p}_1$};
        \node[inputNode, thick] (p3) at (4, 0.25) {$\tilde{p}_2$};
        \node (dots) (p4) at (4, -0.25) {$\vdots$};
        \node[inputNode, thick] (p5) at (4, -1) {$\tilde{p}_d$};
        \node[inputNode, thick] (p6) at (4, -1.75) {$\tilde{q}_1$};

        \draw[stateTransition] (2.4,1) -- (p1);
        \draw[stateTransition] (2.4,0) -- (3.8,0);
        \draw[stateTransition] (2.4,-1) -- (p6);

        \draw (3.6,1.3) rectangle (4.4,-1.3);
        \draw (5.6,1.3) rectangle (6.4,-1.3);
        \draw (5.6,1.6) rectangle (6.4, 2);
        \draw (5.6,-1.6) rectangle (6.4, -2);

        \draw[stateTransition] (p1) -- (5.6,1.8);
        \draw[stateTransition] (p6) -- (5.6,-1.8);
        \draw[stateTransition] (p2) -- (5.6,1);
        \draw[stateTransition] (p3) -- (5.6,0.25);
        \draw[stateTransition] (p5) -- (5.6,-1);
        \draw[stateTransition] (p2) -- (5.6,-1);
        \draw[stateTransition] (p3) -- (5.6,-1);
        \draw[stateTransition] (p2) -- (5.6,0.25);
        \draw[stateTransition] (p5) -- (5.6,0.25);
        \draw[stateTransition] (p3) -- (5.6,1);
        \draw[stateTransition] (p5) -- (5.6,1);

        \node[inputNode, thick] (q1) at (8, 1.8) {$u$};
        \node[inputNode, thick] (q2) at (8, 0.25) {$p_2$};
        \node[inputNode, thick] (q3) at (8, 1) {$p_1$};
        \node (nodes) at (8, -0.25) {$\vdots$};
        \node[inputNode, thick] (q4) at (8, -1) {$p_d$};
        \node[inputNode, thick] (q5) at (8, -1.8) {$q$};

        \draw[stateTransition] (6.4, 1.8) -- (q1);
        \draw[stateTransition] (6.4, -1.8) -- (q5);
        \draw[stateTransition] (6.4, -1) -- (q4);
        \draw[stateTransition] (6.4, 1) -- (q3);
        \draw[stateTransition] (6.4, 0.25) -- (q2);
        \draw[stateTransition] (6.4, 0.25) -- (q3);
        \draw[stateTransition] (6.4, 0.25) -- (q4);
        \draw[stateTransition] (6.4, 1) -- (q2);
        \draw[stateTransition] (6.4, 1) -- (q4);
        \draw[stateTransition] (6.4, -1) -- (q2);
        \draw[stateTransition] (6.4, -1) -- (q3);

        \node[above=of x1, align=center] (l1) {Input \\ Layer};
        \node[right=1.2em of l1, align=center] (l2) {Hidden \\ Layers};
        \node[right=0.5em of l2, align=center] (l3) {Intermediate \\ States};
        \node[right=-0.5em of l3, align=center] (l4) {Decoupled \\Hidden\\ Layers};
        \node[right=0.5em of l4, align=center] (l5) {Output \\ Layer};
    \end{tikzpicture}
\end{center}
\caption{A decoupled neural network structure for solving high-dimensional PDEs with multiple outputs.}
\label{st::nd}
\end{figure}

For high-dimensional problems, a decoupled neural network can better distinguish derivatives and the solution, as is shown in Fig. \ref{st::nd}.
After a certain number of hidden layer operations, the input is transformed into a series of intermediate states. The states are then separated by some independent hidden layers to calculate the different variables. For example, $(p_1,\cdots,p_d)\approx \nabla u$ and $q\approx \Delta u $. The depth of the two types of hidden layers depends on how closely you need the different order derivative to be connected.

\subsection{Local Deep Galerkin Method}
With the system and network structure defined, we can propose local deep learning methods to solve the PDEs with high order derivatives.

Taking the residual of the system (\ref{system}), the modified loss function is given as follows
\begin{equation}
    \begin{array}{rcl}
        J_e(\varphi) &=& \|(\varphi_1)_t - F(\varphi_1,\varphi_2,\cdots,\varphi_{k},D\varphi_{k})\|^2_{\Omega_T} + \sum_{i=1}^{k-1}\|\varphi_{i+1} - D\varphi\|_{\Omega_T}^2,\\
        J_i(\varphi) &=& \|\varphi_1(\cdot,0) - u_0(\cdot)\|^2_{\Omega},\\
    \end{array}
    \label{eq::ddmloss}
\end{equation}
and the boundary term is expressed as
\begin{equation}
    J_b(\varphi)=\left\{
    \begin{array}{ll}
        \|\mathcal{B}\varphi_1 - \mathbf{g}\|^2_{\partial\Omega_T}, &\text{Dirichlet boundary condition,}\\
        \|\mathcal{B}(\varphi_2,\cdots,\varphi_k)^T - \mathbf{g}\|^2_{\partial\Omega_T}, &\text{Neumann boundary condition,}\\
        \|\varphi(\mathbf{x}_l,\cdot) - \varphi(\mathbf{x}_r,\cdot)\|^2_{\partial\Omega_T},&\text{Periodic boundary condition,}
    \end{array}
    \right.
    \label{eq::ddmloss2}
\end{equation}
The local deep Galerkin method, which follows deep Galerkin method\cite{Sirignano2017DGM}, is summarized in Algorithm \ref{al::ddm}. This method treats the solution and its derivatives or other necessary intermediate variables as unknown functions while simultaneously learning their values. These restrictions cause a certain increase in calculations, but this is still far less expansive than calculating the derivatives.

\begin{algorithm}
\caption{Local Deep Galerkin Method (LDGM)}
\begin{algorithmic}[1]
    \State Build up the neural network $\vec{\varphi}(x,t;\theta)$. Determine the hidden layer structure $\mathcal{N}$, the layer width $n$, the dimension of output $m$ and the active functions $\sigma$ according to the given PDEs.
    \State Initialize the parameters $\theta=\theta^0$, sampling times $s_1$, the number of samples $N_e,N_i,N_b$ in $\Omega_T,\Omega,\partial\Omega_T$, optimization steps $s_2$ and the learning rate $\gamma$.
    \For{$i=0:s_1$}
    \State Obtain random sampling points $\{(x_e,t_e)\}_{N_e},\{(x_i,0)\}_{N_i},\{(x_b,t_b)\}_{N_b}$.
    \State Set $\theta^{i,0}=\theta^i$.
    \For{$j=0:s_2$}
        \State Calculate the loss function $J(\varphi(x,t;\theta^{i,j}))$ at sampling points according to (\ref{eq::ddmloss}) and (\ref{eq::ddmloss2}).
        \State Optimize the parameters $\theta$
        \begin{equation*}
            \theta^{i,j+1}=\theta^{i,j} - \gamma \nabla_\theta J(\theta^{i,j}).
        \end{equation*}
    \EndFor
    \State Set $\theta^{i+1} = \theta^{i,s_2+1}$.
    \EndFor
    \end{algorithmic}
    \label{al::ddm}
\end{algorithm}

\subsection{Local Deep Ritz Method}
We can also combine the technique discussed above with the deep Ritz method\cite{weinan2017deep}. Consider the following bi-Laplacian equation:
\begin{equation}
    \left\{
    \begin{aligned}
        \Delta^2 u&= f,&& x\in \Omega,\\
        u(x) &= 0, &&x \in \partial\Omega,\\
        \frac{\partial u}{\partial \mathbf{n}}&=0,&&x\in\partial \Omega,
    \end{aligned}
    \right.
    \label{eq::variational}
\end{equation}
the weak formulation of which is
\begin{equation}
    J(u)=\int_\Omega\left(\frac{1}{2}(\Delta u(x))^2 - f(x)u(x)\right)dx,
\end{equation}
where $u\in H$ and $H$ is the set of trial functions. Using a neural network $\varphi$, the loss function in the deep Ritz method is defined as
\begin{equation}
    J(\varphi) = \int_\Omega \left(\frac{1}{2}\left(\Delta_x\varphi(x;\theta)\right)^2 - f(x)\varphi(x;\theta)\right)dx + \lambda\int_{\partial\Omega}\left(
    \varphi(x;\theta)^2+(\frac{\partial\varphi(x;\theta)}{\partial \mathbf{n}})^2\right)dx.
\end{equation}

We use a multi-output neural network $\varphi(x;\theta):\mathbb{R}^d\times\Theta^M\mapsto\mathbb{R}^{d+1}$, and reformulate the loss function as follows
\begin{equation}
    \begin{aligned}
        \hat{J}(\varphi)&= \int_\Omega \left(\frac{1}{2}(\nabla_x \cdot \mathbf{q}(x;\theta))^2 - f(x)p(x;\theta)+\|\nabla_x p(x;\theta) - \mathbf{q}(x;\theta)\|^2\right)dx \\
        &+\lambda\int_{\partial\Omega} (p(x;\theta)^2+(\mathbf{q}(x;\theta)\cdot\mathbf{n})^2)dx,
    \end{aligned}
\end{equation}
where $p=\varphi_1, \mathbf{q} = (\varphi_2,\cdots,\varphi_{d+1})$. A local deep Ritz method can then be used to solve the variational problem.

Notice that we can solve a fourth order problem with a $d+1$-dimensional output neural network. For high-dimensional problems which are the kind of variational problems, the local deep learning method can also be applied.

\section{Advantages of LDLM}\label{sec4}
In this section, we illustrate the advantages of using the local deep learning method to solve differential equations with high order derivatives.
The intuitive performance of some numerical tests will be shown in Section \ref{sec5}.

\subsection{Reduction of Calculations}
The computational complexity of computing a $k$-th order derivative is about $O(L^kn^{2k})$, which is delineated in Section \ref{sec::complexity}. In the LDLM, a $k$-th order derivative becomes $k-1$ restrictions and $k$ first order derivatives. As the restrictions are much cheaper than computing derivatives, the total cost is about $O(Lkn^2)$. The linear growth with respect to the order $k$ is especially suitable for solving high order PDEs.

\subsection{Improving Robustness for Training Process}
The robustness of local deep learning methods is made manifest in the solving of complex nonlinear differential equations, which contain different order derivatives and multiple scales, like the Cahn--Hilliard equation (\ref{eq::CH}). In Section \ref{sec::complexity}, we see that the scale of a $k$-th order derivative is about
\begin{equation}
    u^{(k)}(x;\theta) \sim (\sigma^{(k)}(x)\theta^k)^L.
    \label{eq::scale}
\end{equation}
Considering the $p$-th order derivative and the $q$-th order derivative in one equation, the formulation contains
\begin{equation}
    \theta^{pL} ((\sigma^{(p)}(x))^L+(\sigma^{(q)}(x))^L\theta^{(q-p)L}).
    \label{eq::scale2}
\end{equation}
Assuming $\sigma^{(p)}(x)\sim \sigma^{(q)}(x)$, the term $\theta^{(q-p)L}$ is sensitive if $p\neq q$ and $L\gg1$. The initialization of the parameters $\theta$ will seriously impact the performance of a neural network. This will cause the information of one of the derivatives to be neglected during the training process when $|(q-p)L|\gg1$.

Additionally, high order derivative terms are often accompanied by small coefficients, like the viscosity $\mu$ in the Navier-Stokes equation and the interface width $\epsilon$ in the Cahn--Hilliard equation. In traditional numerical methods, high order numerical schemes or multi-scale analysis methods can overcome the equation's parameter sensitivity.
But small coefficients, coupled with the different order derivatives of the neural network, can cause great difficulties in optimization. In other words, the neural network may not converge due to the sensitivity of the small coefficients.

The modified loss function of the system only includes
\begin{equation}
    \sigma^L(x)+(\sigma'(x))^L\theta^{L}.
\end{equation}
It is easier to assume $\sigma(x) \sim \sigma'(x)$, and the parameter scale can be balanced by $\sigma(x)\sim \sigma'(x)\theta$. In training, different order derivatives only affect adjacent ones, which leads to a more robust result. 

\subsection{Weaker Restrictions on Active Functions}
Following the above, local deep learning methods place less of restriction on active functions. Different active functions have different uses in neural networks. For example, the ReLU active function allows us to circumvent the gradient vanishing problem and the hyperbolic tangent active function can provide smoothness. Choosing a suitable active function for a given task is an open hyperparameter learning problem.

For PDEs containing the $k$-th order derivative, the neural network $\varphi$ should, at least, belong to $C^{k+1}(\Omega)$ from (\ref{eq::scale}), where $C(\Omega)$ is the collection of continuous functions on $\Omega$ and $C^k(\Omega):=\{ f|f^{(k)}\in C(\Omega)\}$.
Combining this with (\ref{eq::scale2}), the active function should satisfy
\begin{itemize}
    \item[i.]   $ \sigma(x)\in C^{k+1}(\Omega)$;
    \item[ii.]  $\sigma^{(p)}(x)\sim \sigma^{(q)}(x), \forall x\in\Omega, 1\leq p,q\leq k+1. $
\end{itemize}
Condition (i) can be satisfied with smooth nonpolynomial active functions, like the hyperbolic tangent and sigmoid. But a common problem is that these active functions inevitably lead to gradient vanishing and exploding problems. Other popular active functions, like ReLU and ReCU, do not meet the condition. Condition (ii) is much more stringent. Among the common elementary functions, only the exponential function satisfies this condition, but it is not usually used as an active function.

For local deep learning methods, these conditions are weakened as
\begin{itemize}
    \item [iii.] $\sigma(x)\in C^{2}(\Omega)$;
    \item [iv.] $\sigma(x)\sim \sigma'(x), \forall x\in\Omega$.
\end{itemize}
If we approximate the weak solution of the equation, Condition (iii) can be further weakened as $\sigma\in C^1(\Omega)$ and the weak derivative of $\sigma'$ exists. Then, some active functions with nonexistent second order derivatives can be used to train the neural network. Condition (iv) means that $\sigma'$ should be bounded in $\Omega$.

The weaker restriction on the active functions provides more choices in local deep learning methods.

\subsection{More Flexible Choices of Network Structures}
Deep neural networks usually have a large number of hidden layers. As mentioned above, the high order derivatives of deep neural networks not only increase the size of the calculation, but also destabilize the training process.
So local deep learning methods, which have successfully avoided the calculation of high order derivatives, can be better combined with deep neural networks.

In addition, we can also use more complex neural network models. For example, one advantage of the residual layer
\begin{equation}
    \mathcal{N}(x) = \sigma(W_2\sigma(W_1x+b_1)+b_2) + x,
\end{equation}
is that it avoids the gradient vanishing problem. This is because the linear term $x$ keeps an additive constant gradient in the first order derivatives, i.e.,
\begin{equation}
    \frac{\partial \mathcal{N}}{\partial z} = (F'(x) + 1)\frac{\partial x}{\partial z},
\end{equation}
where $F(x) = \sigma(W_2\sigma(W_1x+b_1)+b_2)$. But for high order derivatives, it is not usually effective. For example,
\begin{equation}
    \frac{\partial^3 \mathcal{N}}{\partial z^3} = F'''(x)(\frac{\partial x}{\partial z})^3 + 3F''(x)\frac{\partial x}{\partial z}\frac{\partial^2 x}{\partial z^2} + F'(x)\frac{\partial^3 x}{\partial z^3} + \frac{\partial^3 x}{\partial z^3}.
\end{equation}
it can always avoid the gradient vanishing problem, but it cannot grasp the contribution of various components to the high order derivatives well, which, in turn, will lead to an inaccurate calculation of the derivatives.
For other complex network structures, like the long short term memory layer, the calculation of high order derivatives relies on a series of complex composite functions which greatly increases complexity of the caomputation. For local deep learning methods, which are similar to classical deep learning problems, many existing tools, methods and network structures can be directly migrated to solve high order differential equations.

\section{Numerical Examples}\label{sec5}

\subsection{Setup}
 In this section, we use the local deep Galerkin method (LDGM) to compute a series of examples including high-dimensional linear and nonlinear differential equations with high order derivatives. The accuracy of the solution $\varphi(x;\theta)$ is measured by the relative $L_2$ error $\frac{\|\varphi-u^*\|_2}{\|u^*\|_2}$ where $u^*$ is the exact solution and $\|u\|_2^2=\int_\Omega u^2dx$. If the exact solution is not analytical, we will compare the solution to the reference solution obtained by the finite difference method. The base solution for comparison is obtained by the deep Galerkin method (DGM)\cite{Sirignano2017DGM}.

 The numerical implementation of the algorithm is based on TensorFlow, which is a widely-used open-source software library in machine learning\cite{abadi2016tensorflow}. The automatic differentiation is included in function $\text{tf.gradient}(y,x)$, which returns $\sum_{i=1}^m(\frac{y_i}{x_1},\cdots,\frac{y_i}{x_d})$, rather than the Jacobi matrix.
 In all numerical experiments, a fully connected feedforward network is chosen as the network structure. Unless otherwise noted, the neural network is configured to have $3$ hidden layers with $50$ neurons per hidden layer. Parameters are initialized by the Xavier Initializer (also known as the Glorot Uniform Initializer), which is used to avoid the gradient vanishing and exploding problems\cite{glorot2010understanding}. Most active functions are selected as $\tanh(x)$ for smoothness and final active functions are determined by practical problems. In optimization, we set the learning rate $r = 0.001$, sampling times $s_1 = 1000$, optimization steps $s_2 = 5$ as a default and use the Adam optimizer. Sampling settings are $N_e = 200, N_i = 50, N_b=50$, which contain a total of 300 nodes per suboptimization problem, and the uniform distribution is used for sampling. The weights of loss functions are chosen equal $w_e=w_i=w_b=1$, unless there is a singularity on the boundary or initial condition.

 Notations of the experiments and algorithm parameters are summarized in Table (\ref{tb::notations}) for quick reference.
 \begin{table}[htbp]
    \centering
	\caption{The List of Parameters}
    \begin{tabular}{|c|l|}
    \hline
    Notation & Stands for ... \\
    \hline
    $\varphi(x,t;\theta)$ & Neural network of input $(x,t)$ with trainable parameters $\theta$\\
    \hline
    $d$ & Dimension of $\Omega \subset \mathbb{R}^d$\\
    \hline
    $L$ & Number of hidden layers \\
    \hline
    $m$ & Dimension of output \\
    \hline
    $n$ & Hidden layer width \\
    \hline
    $\sigma$ & Active functions in neural network\\
    \hline
    $J_e,J_i,J_b$ & The loss of the equation, the initial value and the boundary condition\\
    \hline
    $w_e,w_i,w_b$& The weights of the losses\\
    \hline
    $N_e,N_i,N_b$& Number of sampling nodes on the domain $\Omega_T$,$\Omega\times\{0\}$ and $\partial\Omega_T$\\
    \hline
    $r$ & Learning rate of network parameter $\theta$ \\
    \hline
    $s_1$ & Sampling times on the whole training process\\
    \hline
    $s_2$ & Optimization steps of per sampling stage\\
    \hline
    \end{tabular}
    \label{tb::notations}
 \end{table}

\subsection{Fourth Order PDE}
In the first case, we show that while there is no obvious difference between the LDGM and the DGM in terms of accuracy, the LDGM can greatly speed up the calculation.

We consider a simple model for a vibrating elastic beam first\cite{niordson1965optimal}:
\begin{equation}
    u_t = -u_{xxxx},\quad x\in[0,2\pi],\quad t\in[0,1].
\end{equation}
With the Dirichlet boundary conditions $u(x,t) = 0,u_{xx}(x,t) = 0, x\in\partial\Omega$ and the initial condition $u_0(x)=\sin x$, the exact solution is given as $u(x,t) = e^{-t}\sin x$. It is costly to calculate the fourth order derivative in deep neural networks while we get high order derivatives directly from the multi-output neural network.

\begin{figure}[htbp]
    \centering
    \includegraphics[width=\textwidth]{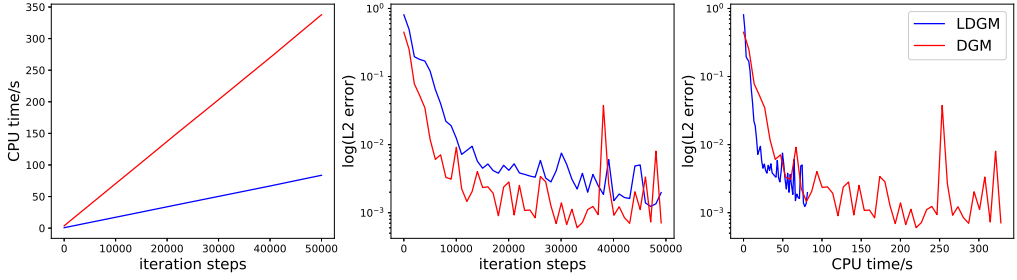}\vspace{-1em}
    \caption{The training processes of DGM and LDGM with the learning rate $r=10^{-4}$ and $s_1\times s_2=50000$ training steps. The red line is the deep Galerkin method and the blue line is the local deep Galerkin method. The left figure shows the time spent between two methods under the same iteration step. The middle figure shows the logarithmic $L_2$ error with respect to the iteration step. The right figure shows how fast the $L_2$ error drops.}
    \label{fig::fourth-PDE}
\end{figure}

From Fig. \ref{fig::fourth-PDE}, we find that under the same iteration step, the LDGM is trained much faster than the DGM, while the error is slightly different. In the third picture, we can conclude that the LDGM approaches the solution of the PDEs containing high order derivatives faster than the DGM. It saves a lot of time in training, and the advantage will be magnified as the network's depth increases.

Another interesting aspect of this study is that the oscillation amplitude of the LDGM is smaller than that of the DGM. Although not sufficient, we realize that the local deep Galerkin method is more stable
for PDEs with high order derivatives.

\subsection{Cahn--Hilliard Equation}
\label{sec::CH}
The Cahn--Hilliard(CH) equation is a popular mathematical physical equation used to describe the process of phase separation. When we use the deep learning method to solve the Cahn--Hilliard equation, it fails when $\epsilon$ is small. This is why we propose the local method to strengthen the robustness in training.

\begin{figure}[!t]
    \centering
    \subfigure[$\epsilon=0.1$]{\includegraphics[width=\textwidth]{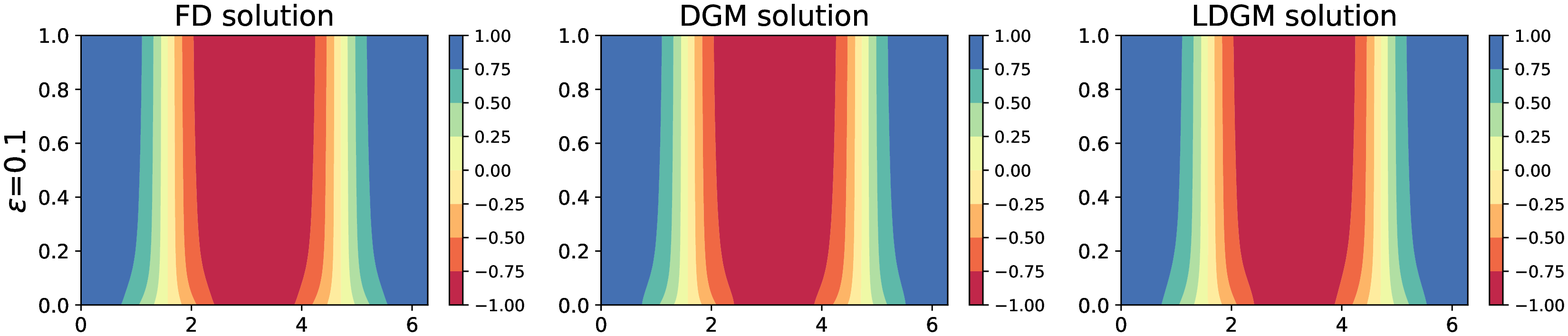}}
    \subfigure[$\epsilon=0.03$]{\includegraphics[width=\textwidth]{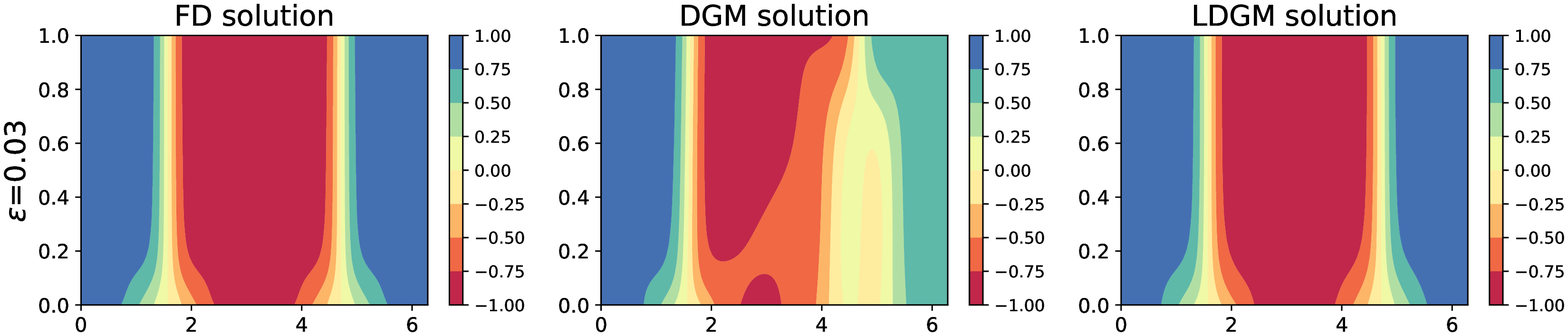}}
    \subfigure[$\epsilon=0.01$]{\includegraphics[width=\textwidth]{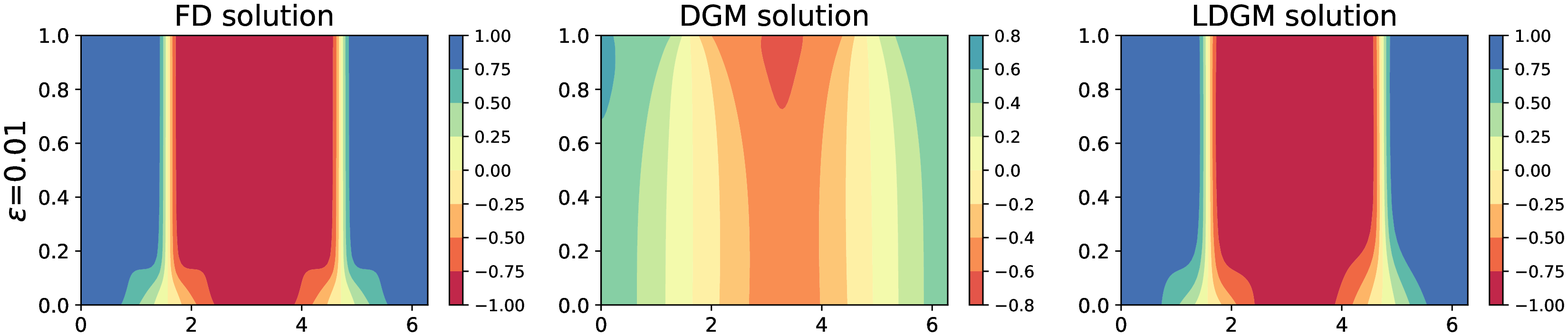}}
    \caption{The solutions to the Cahn--Hilliard equation are by the finite difference method, the DGM and the LDGM from left to right with different $\epsilon=0.1,0.03,0.01$. When $\epsilon$ is small, the DGM fails to approximate the solution under given sampling stages $s_1 =5000$.}
    \label{fig::CH}
\end{figure}

The 1-D CH equation can be given as
\begin{equation}
    u_t + \epsilon u_{xxxx} + f(u)_{xx} = 0, \quad x\in [0,2\pi], \quad t\in[0,1],
    \label{eq::CH}
\end{equation}
where $f(u)=u-u^3$. Given the initial condition $u_0(x)=\cos x$ and the zero Neumann boundary condition, we define the following loss function
\begin{equation}
    \begin{array}{rcl}
        J_e(\varphi) &=& \|(\varphi_1)_t - (\varphi_4)_x\|_{\Omega_T}^2 + \|\varphi_3 +\epsilon (\varphi_2)_{x} + f(\varphi_1)\|_{\Omega_T}^2 + \|\varphi_2 - (\varphi_1)_x\|_{\Omega_T}^2 + \|\varphi_4-(\varphi_3)_x\|_{\Omega_T}^2,\\
        J_i(\varphi) &=& \|\varphi_1 - \cos(x)\|_{\Omega}^2,\\
        J_b(\varphi)&=& \|\varphi_2\|_{\partial\Omega_T}^2 + \|\varphi_4\|_{\partial\Omega_T}^2,\\
        J(\varphi)&=& J_e(\varphi) + J_i(\varphi) + J_b(\varphi).
    \end{array}
\end{equation}
Here, the multi-output neural network is
\begin{equation}
    \varphi(x,t;\theta) = \begin{pmatrix}
        \varphi_1\\
        \varphi_2\\
        \varphi_3\\
        \varphi_4
    \end{pmatrix} \approx \begin{pmatrix}
        u\\
        u_x\\
        \phi\\
        \phi_x
    \end{pmatrix},
\end{equation}
where $\phi = -\epsilon u_{xx} - f(u)$. For the DGM, classical $L_2$ loss is used. The reference solution is given by a finite difference method with 129 spectral nodes, and the numerical scheme is
\begin{equation}
    \frac{\hat{u}_k^{n+1}-\hat{u}_k^n}{\delta t} + \epsilon (ik)^4 \hat{u}_k^{n+1} + (ik)^2\hat{f}^n_k=0,
\end{equation}
where $f^n = u^n - (u^n)^3$, $\delta t =0.01$. The spectral method and the FFT solver are used here. We show numerical results in Fig. \ref{fig::CH} with sampling stages $s_1= 5000$.

%\vspace{-1em}
%\newpage
%\vspace{-2em}

When $\epsilon \sim O(1)$ or there is a source term $g(x,t,\epsilon)$, the DGM can be guided to approximate the true solution quickly.
But when $\epsilon$ becomes small, the DGM fails while the LDGM is still able to capture the interface. For the DGM, it costs about $s_1 = 40000$ sampling stages to get a reasonably accurate solution when $\epsilon=0.02$ and it costs about $s_1 =100000$ sampling stages when $\epsilon=0.01$ in experiments.

The following reasons together result in this failure:
\begin{itemize}
    \item[a.] The Xavier initializer gives $\text{Var}(\theta) = \frac{1}{n}$ and $E(\theta) = 0$. So the fourth order term $\epsilon u_{xxxx}$ modeled by $\epsilon(\theta^4\sigma^{4})^L$ leads to a gradient vanishing problem at the beginning.
    \item[b.] The number of sampling nodes is insufficient to capture the interface. In each suboptimization problem, only about $\frac{\epsilon N_e}{2\pi}$ nodes are around the interface. It follows that the parameters are updated slowly.
    \item[c.] It is a non-convex optimization, and the learning rate needs to be small, which causes the parameters to stay in an incorrect interval for a long time.
\end{itemize}

A specific initialization and sampling method may work, but it is not a general strategy for different PDEs. In the LDGM, the order of derivatives has been reduced, which makes this algorithm less affected by the gradient vanishing problem caused by initialization. In addition, changing $\epsilon$ will not significantly affect the convergence speed of the LDGM. In other words, the robustness of the LDGM is better than that of the DGM.

\subsection{Modified KdV Equation}
In this test, we show that when the neural network becomes deeper, the parameter scale difference between different order derivatives becomes more obvious. Consider the following modified Korteweg-de Vries equation
\begin{equation}
    \left\{\begin{aligned}
        &u_t -6u^2u_x+u_{xxx}=0, x\in[-2,2], &&t\in[0,1],\\
        &u(2,t) = \tanh(2t+1),u(-2,t)=\tanh(2t-3), &&t\in[0,1],\\
        &u(x,0) = \tanh(x-1), &&x\in[-2,2].
    \end{aligned}\right.
    \label{eq::kdv}
\end{equation}

The kink solution of problem (\ref{eq::kdv}) is $u(x,t) = \tanh(x+2t-1)$.
Set up the neural network
\begin{equation}
    \varphi(x,t;\theta) = \begin{pmatrix}
        \varphi_1\\ \varphi_2 \\ \varphi_3
    \end{pmatrix}\approx
    \begin{pmatrix}
        u\\ u_x \\ u_{xx}
    \end{pmatrix},
\end{equation}
and the loss function is given as follows
\begin{equation}
    \begin{array}{rcl}
        J_e(\varphi) &=& \|(\varphi_1)_t - 6\varphi^2_{1}\varphi_2 + (\varphi_3)_x\|_{\Omega_T}^2 + \sum_{i=1}^2\|(\varphi_i)_x-\varphi_{i+1}\|_{\Omega_T}^2,\\
        J_i(\varphi) &=& \|\varphi_1 - \tanh(x-1)\|_{\Omega}^2,\\
        J_b(\varphi)&=& \|\varphi_1(-2,t) - \tanh(2t-3)\|_{[0,1]}^2+\|\varphi_1(2,t) - \tanh(2t+1)\|_{[0,1]}^2.
    \end{array}
\end{equation}
\begin{figure}[htbp]
    \centering
    \includegraphics[width=0.9\textwidth]{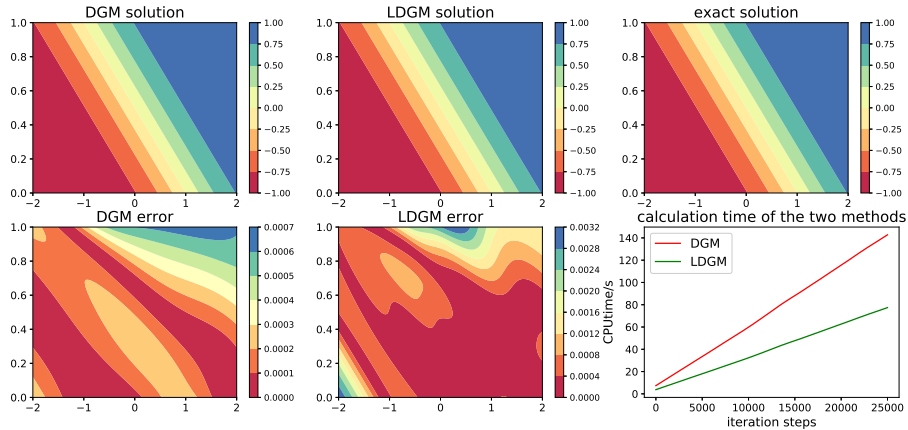}\vspace{-1em}
    \caption{Deep learning solutions of the modified KdV equation obtained by the DGM and the LDGM. Under the default settings, the LDGM is faster than the DGM. }
    \label{fig::kdv}
\end{figure}

Under the default settings in Fig. \ref{fig::kdv}, the only proven superiority of the LDGM is that it costs less time when solving the problem. We then compare the performance of different network settings in solving the modified KdV equation. The relative $L_2$ errors after 25000 steps training are shown in Table \ref{tb::kdv}.

Table \ref{tb::kdv} shows that the LDGM is less affected by network settings, while the DGM is more sensitive. When $L=48, n= 10$, the DGM fails to solve the KdV equation due to the gradient vanishing problem caused by the active function $\tanh(x)$. But under the continuity assumption of the DGM, the active function of the neural network does not have many choices. In the LDGM, it only requires $\sigma(x)\in \mathbb{C}^2$. Here the exponential linear units
\begin{equation}
    \sigma(x) = \left\{\begin{array}{ll}
        x, &x>0,\\
        \alpha (e^x-1), &\text{otherwise}.
    \end{array}\right.
\end{equation}
can be used to further reduce errors and avoid gradient vanishing problem. When $L=64$, we repeat the experiment 100 times and record whether the relative $L_2$ error is less than $1\%$ after 10000 iteration steps. The success rates of the DGM and the LDGM are $30\%$ and $89\%$, respectively. This means that when the network is deep, the DGM is almost ineffective, while the LDGM provides more flexible choices in terms of network structures and active functions.
\begin{table}[htbp]
    \centering
	\caption{Relative $L_2$ Errors under Various Network Coefficients}
    \begin{tabular}{|c|c|c|c|c|c|c|c|}
    \hline
    $n=10$,various $L$ & 3  & 6 & 9 & 12 & 24 &48 \\
    \hline
    DGM & 0.132\% & 0.099\% & 0.155\% &0.171\% &0.579\% & 84.22\% \\
    \hline
    LDGM & 0.196\% & 0.164\% &0.165\%& 0.142\% & 0.098\% & 0.243\%  \\
    \hline
    $L=3$,various $n$ & 10 & 20 & 40 & 80 & 160 & 320 \\
    \hline
    DGM & 0.132\% & 0.059\% & 0.060\% &0.209\% &0.027\% & 1.694\% \\
    \hline
    LDGM & 0.196\% & 0.304\% &0.174\%& 0.210\% & 0.191\% & 0.342\%  \\
    \hline
    \end{tabular}
    \label{tb::kdv}
\end{table}

\subsection{High-Dimensional Heat Equation}
Innumeriable previous studies have shown that deep learning methods have distinct advantages in solving high-dimensional problems. In this example, we show that the LDGM inherits these advantages.
Consider the general heat equation
\begin{equation}
    \begin{array}{rcll}
        u_t - \Delta u &=& f(x,t), &x \in \Omega=[0, 1]^d, t\in [0,1],\\
        u(x,t) &=& g(x,t), &x \in \partial\Omega, t\in[0,1],
    \end{array}
    \label{eq::heat-nd}
\end{equation}
where $f(x,t) = 2d(t+1) + \sum_{i=1}^d x_i(1-x_i)$ and $g(x,t) = \sum_{i=1}^d x_i(1-x_i)(t+1) $. Given the initial condition $u_0(x) = \sum_{i=1}^d x_i(1-x_i)$, problem (\ref{eq::heat-nd}) has a classical solution
$u^*(x,t) = \sum_{i=1}^d x_i(1-x_i)(t+1)$. Define the multi-output neural network as follows
\begin{equation}
    \varphi(x,t;\theta) = \begin{pmatrix}
        \varphi_1\\
        \varphi_2\\
        \vdots\\
        \varphi_{d+1}
    \end{pmatrix}\approx\begin{pmatrix}
        u\\
        u_{x_1}\\
        \vdots\\
        u_{x_d}
    \end{pmatrix}.
\end{equation}
The loss function has the form of
\begin{equation}
    \begin{array}{rcl}
        J_e(\varphi) &=& \|(\varphi_1)_t - \sum_{i=1}^d (\varphi_{i+1})_{x_i} -f\|_{\Omega_T}^2 + \sum_{i=1}^d\|(\varphi_1)_{x_i}-\varphi_{i+1}\|_{\Omega_T}^2,\\
        J_i(\varphi) &=& \|\varphi_1 - u_0\|_{\Omega}^2,\\
        J_b(\varphi)&=& \|\varphi_1 - g\|_{\partial\Omega_T}^2.
    \end{array}
\end{equation}
Notice that we only need to output the first order derivatives of all dimensions for computing the second order derivatives, which greatly saves calculation and storage space. The same strategy applies to higher order derivatives.

\begin{figure}[!t]
    \centering
    \includegraphics[width=\textwidth]{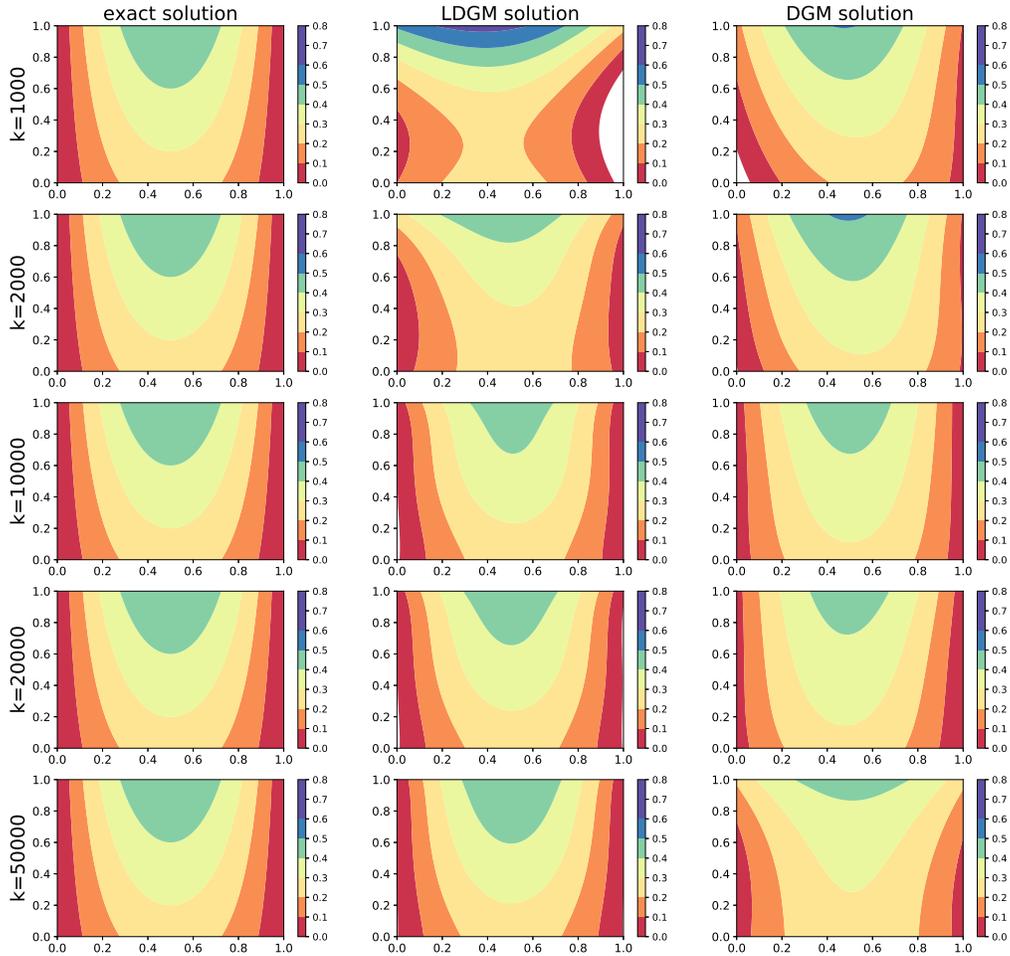}\vspace{-5em}
    \caption{The training process of the high-dimensional heat equation obtained by DGM and LDGM with respect to dimensions $d=5$. The picture shows the exact solution, the solution of DGM and the solution of LDGM from left to right. From top to bottom, it shows the solutions when the iteration steps $k=1000,2000,10000,20000,50000$. For display purposes, the images show the slices of $x_2,x_3,\cdots,x_d=0$. The abscissa is $x_1$ and the ordinate is $t$.}
    \label{fig::ndheat}
\end{figure}

\begin{figure}[htbp]
    \centering
    \includegraphics[width=\textwidth]{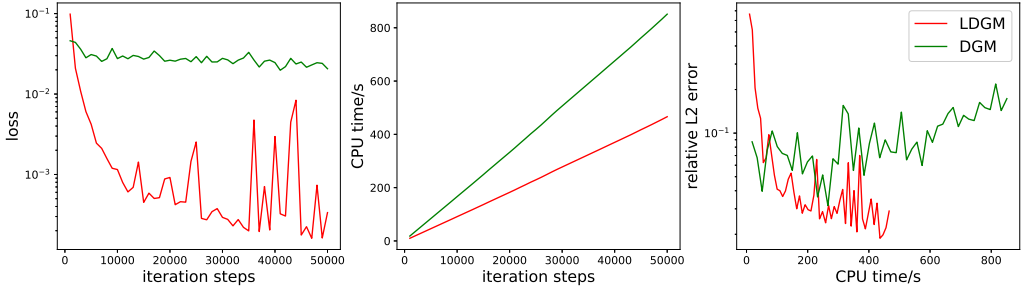}
    \caption{The iteration curves of the DGM and the LDGM for the 5-D heat equation. From left to right: loss vs. iteration steps, relative $L_2$ error vs. iteration steps, and relative $L_2$ error vs. calculation time, respectively.}
    \label{fig::heat_loss}
\end{figure}

Setting $d = 5$, $L=4$, $n=100$, $r=0.0005$  and using 50000 iteration steps ($s_1 = 10000$ and $s_2=5$), the solutions are given in Fig. \ref{fig::ndheat}.
%{-1em}
Although the LDGM needs to learn more details in optimization, the DGM has no convergence after a reasonable iteration step. The error of the LDGM is caused by the limited approximation ability of such a neural network and training set. The default settings is not sufficient to cover the entire region $\Omega^d$. For a more precise solution, adding hidden layers, expanding the network's width and increasing the number of sampling nodes are all viable options. Fig. \ref{fig::heat_loss} uses the error curve to provide a more intuitive comparison. High-dimensional second order derivatives calculated by the automatic differentiation in the loss are hard to optimize. Sometimes, an adaptive piecewise learning rate like $r \sim 10^{-[log(k)]}$ is chosen, where $k$ is the iteration step. It always works but it can be costly. In the LDGM, the error drops quickly, which means that the LDGM retains its advantages in solving high-dimensional problems.

\section{Concluding Remarks \& Declaration}
In this paper, we list the difficulties associated with computing high order derivatives of neural networks for solving PDEs. Calculating high order derivatives is costly and can cause a parameter scaling problem. In addition, complex calculations limit our choices of network structures and active functions. We propose a local deep learning method to overcome these problems.
 We consider the derivatives of the solution as intermediate variables and rewrite the original problem as a system of low order PDEs. The loss function takes the residual of the equivalent system.
 With a multi-output neural network, the local deep learning method is established. We demonstrate the performance of the local deep Galerkin method on a variety of PDEs, including high-dimensional problems, phase field problems and high order PDEs. In all numerical tests, the local deep Galerkin method is shown to be both stable and highly efficient.
 
 {\bf Declaration:} We enclose this paper by a declaration to show the originality of our work. Our project started more than one and half a years ago. About finishing the work by the end of 2020, we found the work \cite{lyu2020mim} on arxiv post in June 2020. The exactly same technique is adopt to solve high-order PDEs by deep learning. They called it as deep mixed residual method, in contrast called as local deep learning method in our paper. Next we try to show the independence and originality of our work by shortly presenting how we propose the so-called local deep leaning method and the difference between our work with \cite{lyu2020mim}. First, our motivation of this study was to use deep learning methods to solve phase field equations. While we successfully solved the Allen--Cahn equation, we failed to solve the Cahn--Hilliard equation, even though they are both gradient flow problems associating with same free energy. We later discovered that our failure to solve the CH equation was due to the high order derivative and small $\epsilon$ causing a parameter scaling problem. Thus the LDLM was proposed. Secondly, \cite{lyu2020mim} focuses more on how the technique results in more accurate solutions as well as more accurate derivatives. But we focus more on how this technique reduces the computations, improves the robustness of the training process, and overcomes the gradient vanishing or exploding problem in high order PDEs.

\section*{Acknowledgements}
This work is partially supported by the National Natural Science Foundation of China/Hong Kong RRC Joint Research Scheme (NSFC/RGC 11961160718), and the fund of the Guangdong Provincial Key Laboratory of Computational Science and Material Design (No. 2019B030301001). The work of J. Yang is supported by the National Science Foundation of China (NSFC-11871264) and the Guangdong Basic and Applied Basic Research Foundation (2018A0303130123).
\bibliographystyle{plain}
\bibliography{ref}

\begin{thebibliography}{10}

\bibitem{abadi2016tensorflow}
Mart{\'\i}n Abadi, Paul Barham, Jianmin Chen, Zhifeng Chen, Andy Davis, Jeffrey
  Dean, Matthieu Devin, Sanjay Ghemawat, Geoffrey Irving, Michael Isard, et~al.
\newblock Tensorflow: A system for large-scale machine learning.
\newblock In {\em 12th $\{$USENIX$\}$ symposium on operating systems design and
  implementation ($\{$OSDI$\}$ 16)}, pages 265--283, 2016.

\bibitem{DBLP:journals/corr/BaydinPR15}
Atilim~Gunes Baydin, Barak~A. Pearlmutter, and Alexey~Andreyevich Radul.
\newblock Automatic differentiation in machine learning: a survey.
\newblock {\em CoRR}, abs/1502.05767, 2015.

\bibitem{bengio1994learning}
Yoshua Bengio, Patrice Simard, and Paolo Frasconi.
\newblock Learning long-term dependencies with gradient descent is difficult.
\newblock {\em IEEE transactions on neural networks}, 5(2):157--166, 1994.

\bibitem{cockburn1998local}
Bernardo Cockburn and Chi-Wang Shu.
\newblock The local discontinuous galerkin method for time-dependent
  convection-diffusion systems.
\newblock {\em SIAM Journal on Numerical Analysis}, 35(6):2440--2463, 1998.

\bibitem{dockhorn2019discussion}
Tim Dockhorn.
\newblock A discussion on solving partial differential equations using neural
  networks.
\newblock {\em arXiv preprint arXiv:1904.07200}, 2019.

\bibitem{DBLP:journals/corr/abs-1710-00211}
Weinan E and Bing Yu.
\newblock The deep ritz method: {A} deep learning-based numerical algorithm for
  solving variational problems.
\newblock {\em CoRR}, abs/1710.00211, 2017.

\bibitem{glorot2010understanding}
Xavier Glorot and Yoshua Bengio.
\newblock Understanding the difficulty of training deep feedforward neural
  networks.
\newblock In {\em Proceedings of the thirteenth international conference on
  artificial intelligence and statistics}, pages 249--256, 2010.

\bibitem{grosse2017lecture}
Roger Grosse.
\newblock Lecture 15: Exploding and vanishing gradients.
\newblock {\em University of Toronto Computer Science}, 2017.

\bibitem{Han2018Solving}
Jiequn Han, Jentzen Arnulf, and E.~Weinan.
\newblock Solving high-dimensional partial differential equations using deep
  learning.
\newblock {\em Proceedings of the National Academy of Sciences}, pages
  201718942--, 2018.

\bibitem{han2020derivative}
Jihun Han, Mihai Nica, and Adam~R Stinchcombe.
\newblock A derivative-free method for solving elliptic partial differential
  equations with deep neural networks.
\newblock {\em arXiv preprint arXiv:2001.06145}, 2020.

\bibitem{hanin2018neural}
Boris Hanin.
\newblock Which neural net architectures give rise to exploding and vanishing
  gradients?
\newblock In {\em Advances in Neural Information Processing Systems}, pages
  582--591, 2018.

\bibitem{hayati2007feedforward}
Mohsen Hayati and Behnam Karami.
\newblock Feedforward neural network for solving partial differential
  equations.
\newblock {\em Journal of Applied Sciences}, 7(19):2812--2817, 2007.

\bibitem{hochreiter1998vanishing}
Sepp Hochreiter.
\newblock The vanishing gradient problem during learning recurrent neural nets
  and problem solutions.
\newblock {\em International Journal of Uncertainty, Fuzziness and
  Knowledge-Based Systems}, 6(02):107--116, 1998.

\bibitem{Hornik1991Approximation}
Kurt Hornik.
\newblock Approximation capabilities of multilayer feedforward networks.
\newblock {\em Neural Networks}, 4(2):251--257, 1991.

\bibitem{Hornik1989Multilayer}
Kurt Hornik, Maxwell Stinchcombe, and Halbert White.
\newblock Multilayer feedforward networks are universal approximators.
\newblock {\em Neural Networks}, 2(5):359--366, 1989.

\bibitem{kingma2014adam}
Diederik~P Kingma and Jimmy Ba.
\newblock Adam: A method for stochastic optimization.
\newblock {\em arXiv preprint arXiv:1412.6980}, 2014.

\bibitem{lagaris1998artificial}
Isaac~E Lagaris, Aristidis Likas, and Dimitrios~I Fotiadis.
\newblock Artificial neural networks for solving ordinary and partial
  differential equations.
\newblock {\em IEEE transactions on neural networks}, 9(5):987--1000, 1998.

\bibitem{lagaris2000neural}
Isaac~E Lagaris, Aristidis~C Likas, and Dimitris~G Papageorgiou.
\newblock Neural-network methods for boundary value problems with irregular
  boundaries.
\newblock {\em IEEE Transactions on Neural Networks}, 11(5):1041--1049, 2000.

\bibitem{lu2019deepxde}
Lu~Lu, Xuhui Meng, Zhiping Mao, and George~E Karniadakis.
\newblock Deepxde: A deep learning library for solving differential equations.
\newblock {\em arXiv preprint arXiv:1907.04502}, 2019.

\bibitem{lyu2020mim}
Liyao Lyu, Zhen Zhang, Minxin Chen, and Jingrun Chen.
\newblock Mim: A deep mixed residual method for solving high-order partial
  differential equations.
\newblock {\em arXiv preprint arXiv:2006.04146}, 2020.

\bibitem{niordson1965optimal}
Frithiof~I Niordson.
\newblock On the optimal design of a vibrating beam.
\newblock {\em Quarterly of Applied Mathematics}, 23(1):47--53, 1965.

\bibitem{pascanu2013difficulty}
Razvan Pascanu, Tomas Mikolov, and Yoshua Bengio.
\newblock On the difficulty of training recurrent neural networks.
\newblock In {\em International conference on machine learning}, pages
  1310--1318, 2013.

\bibitem{RAISSI2019686}
M.~Raissi, P.~Perdikaris, and G.E. Karniadakis.
\newblock Physics-informed neural networks: A deep learning framework for
  solving forward and inverse problems involving nonlinear partial differential
  equations.
\newblock {\em Journal of Computational Physics}, 378:686 -- 707, 2019.

\bibitem{Sirignano2017DGM}
Justin Sirignano and Konstantinos Spiliopoulos.
\newblock Dgm: A deep learning algorithm for solving partial differential
  equations.
\newblock {\em Journal of Computational Physics}, 375, 2017.

\bibitem{weinan2017deep}
E~Weinan, Jiequn Han, and Arnulf Jentzen.
\newblock Deep learning-based numerical methods for high-dimensional parabolic
  partial differential equations and backward stochastic differential
  equations.
\newblock {\em Communications in Mathematics and Statistics}, 5(4):349--380,
  2017.

\bibitem{xu2010local}
Yan Xu and Chi-Wang Shu.
\newblock Local discontinuous galerkin methods for high-order time-dependent
  partial differential equations.
\newblock {\em Communications in Computational Physics}, 7(1):1, 2010.

\bibitem{yan2002local}
Jue Yan and Chi-Wang Shu.
\newblock A local discontinuous galerkin method for kdv type equations.
\newblock {\em SIAM Journal on Numerical Analysis}, 40(2):769--791, 2002.

\bibitem{yan2002ldg}
Jue Yan and Chi-Wang Shu.
\newblock Local discontinuous galerkin methods for partial differential
  equations with higher order derivatives.
\newblock {\em Journal of Scientific Computing}, 17(1-4):27--47, 2002.

\bibitem{ZANG2020109409}
Yaohua Zang, Gang Bao, Xiaojing Ye, and Haomin Zhou.
\newblock Weak adversarial networks for high-dimensional partial differential
  equations.
\newblock {\em Journal of Computational Physics}, 411:109409, 2020.

\end{thebibliography}

%% Cover Letter
%In this work, we propose a local deep learning method for solving high order differential equations. The technique is simply to reduce the order of derivatives by introducing intermediate variables. Thus, the high order differential equations can be written into a system low order differential equations. This method is inspired by the idea of local discontinuous Galerkin method. So we call it as local deep learning method. Then we solve the resulted system by a multi-output deep neural network. Direct advantages of this method are reduction of computations, weaker restrictions on active functions and, more flexible choices of network structures, due to avoiding computing high order derivatives of the network. Besides, numerous numerical examples demonstrate how this technique improves the robustness of the training process and overcomes the gradient vanishing or exploding problem in high order PDEs. We also try to provide some brief insights and explanations to these advantages.
\end{document}